\documentclass[final,epsfig,amsfont]{article}
\usepackage{latexsym}
\usepackage[dvips]{graphicx}
\newcommand{\eeq}{\end{equation}}
\newcommand{\beq}{\begin{equation}}
\newcommand{\nuq}[1]{\label{#1} \eeq}
\newcommand{\dbleint}{\int  \!\!\!\!\! \int \! \!}
\newtheorem{definition}{Definition}
\newtheorem{remark}{Remark}
\newtheorem{proposition}{Proposition}
\newtheorem{lemma}{Lemma}
\newtheorem{theorem}{Theorem}
\newtheorem{example}{Example}
\newcommand{\qed}{\rule{2mm}{3mm}}

\begin{document}

\title{Direct and Inverse Computation of Jacobi Matrices of Infinite Homogeneous Affine I.F.S. \thanks{Supported by MIUR-PRIN "Nonlinearity and disorder in classical and quantum transport processes". }
}

\author{
Giorgio Mantica \\
Center for Non-linear and Complex Systems, \\ Department of Physics and Mathematics, \\Universit\`a dell' Insubria, 22100 Como, Italy \\ and CNISM unit\`a di Como,
and  I.N.F.N. sezione di Milano.}

\maketitle

\begin{abstract}
We introduce a new set of algorithms to compute Jacobi matrices associated with measures generated by infinite systems of iterated functions. We demonstrate their relevance in the study of theoretical problems, such as the continuity of these measures and the logarithmic capacity of their support.
Since our approach is based on a reversible transformation between pairs of Jacobi matrices, we also discuss its application to an inverse / approximation problem.
Numerical experiments show that the proposed algorithms are stable and can reliably compute Jacobi matrices of large order. \\

{\em Keywords: Jacobi Matrices -- Orthogonal Polynomials -- Singular Continuous Measures -- Iterated Function Systems -- Gaussian--I.F.S. Quadratures -- Inverse Fractal Problems -- Discrete Schr\"odinger Operators } \\

{\em MSC2000 Class: 42C05 -- 47B36 -- 65D32 -- 65J22 -- 81Q10}
\end{abstract}

\section{Introduction} \label{secintro}
Let $\mu$ be a non-negative probability measure with compact support
in $\bf R$ and let $J_\mu$ be the associated Jacobi matrix:
 \beq
   J_\mu :=
           \left(   \begin{array}{ccccc}  a_0  & b_1 &      &      &       \cr
                               b_1  & a_1 & b_2  &      &       \cr
                                    & b_2 & a_2  & b_3  &       \cr
                  &     &   \ddots    &    \ddots   &   \ddots \cr
                               \end{array} \right) .
 \nuq{jame}
This infinite symmetric tridiagonal matrix encodes the three-terms recurrence relation of
the orthonormal polynomials $\{p_n(\mu;s)\}_{n \in {\mathbf N}}$ of
$\mu$:
\begin{equation}
\label{nor2}
   s p_n (\mu;s) = b_{n+1} p_{n+1}(\mu;s)
   + a_n p_n(\mu;s)  + b_{n} p_{n-1}(\mu;s),
\end{equation}
with $b_0=0$ and initialized by
$p_{-1}(\mu;s) = 0$, $p_0(\mu;s) = 1.$ Recall that the integral $\int p_n(\mu;s) p_m(\mu;s) d \mu(s)$ is equal to one when $n=m$ and is null in all other cases.

According to Gautschi \cite{gaut}, the computation of the Jacobi
matrix associated with a given measure is a {\em fundamental
problem} of numerical analysis. In fact, when the moment problem is determined \cite{ach} (this is the case under the compactness hypothesis above),  the Jacobi matrix $J_\mu$ uniquely identifies the measure $\mu$. Its computation is necessary for the evaluation of orthogonal sums via Clenshaw's algorithm \cite{clen}, but arguably its most important r\^ole comes in differentiation and integration using orthogonal polynomials \cite{tom} and in Gaussian quadratures \cite{gaut2,gaut3,dirk0}, via the linear algebra technique of Golub and Welsch \cite{golw,golmeur}.

Jacobi matrices play a major r\^ole also in mathematical physics, since they can be seen as discrete Schr\"{o}dinger operators acting in $l_2({\bf Z}_+)$, properly when $b_n=1$ for all $n$ and in an extended sense in the general case. The investigation of the links between the properties of the sequences of coefficients $a_n$ and $b_n$ and those of the ``spectral'' measure $\mu$ has involved both students of orthogonal polynomials and numerical analysis  and of mathematical physics and quantum dynamics, see {\em e.g.} \cite{vanas}, \cite{damasi} and references therein. 
In the latter setting the Fourier transform of $\mu$ and of its orthogonal polynomials take on the dynamical meaning of probability amplitudes of the quantum motion \cite{poin1,poin2,etna}. Remark that they can be numerically computed in a very efficient way directly from the Jacobi matrix \cite{ober}.

The asymptotic properties of the Fourier transform of a measure have been the focus of intense investigation also in harmonic analysis \cite{str1,str2,str3,stric}: these studies have highlighted the distinctive properties of singular continuous measures, which appear rarely on the stage of numerical analysis, but that, to the contrary, are principal in this paper.

In fact, we consider herein measures generated by
Systems of Iterated Functions (I.F.S) \cite{hut,dem}, that are frequently, although not always, singular with respect to the Lebesgue measure.
I.F.S. are a highly versatile tool, that has been employed for the
approximation of natural objects as well as for image compression \cite{diaco,ba2}, for wavelet construction \cite{jeff} and for numerical integration \cite{steve}. In this forcefully minimal list of applications of I.F.S. an item of theoretical relevance must be added: Jacobi matrices of I.F.S. measures have been used as abstract models of aperiodic Schr\"odinger operators \cite{gioitalo,physd1,physd2,etna}, whose fine spectral properties are crucial for the phenomenon of wave propagation, in a study which lies
at the intersection of the three disciplines just mentioned: numerical analysis, mathematical physics and harmonic analysis.

Computing the Jacobi matrix associated with an I.F.S. measure is therefore a fundamental task, but a challenging one, the more so because the usual algorithms based on modified moments \cite{sak} are ill conditioned, for reasons explained in \cite{hans1,cap} (see also \cite{berna}). This goal has been achieved with {\em ad hoc} techniques \cite{cap,hans2,mobius,escribano} for I.F.S. with {\em finitely many} maps. In parallel, the computation of the Jacobi matrix for refinable functionals \cite{wsw,gaut2p} (a particular case of I.F.S.)  provided results \cite{dirkvil,dirkvil2} that can compared to those for I.F.S. Yet, not all of the above algorithms are capable of producing Jacobi matrices of large orders, a crucial requirement for the investigation of the theoretical questions mentioned before.

In this paper, we contribute a new set of algorithms to this family, that are numerically stable to large orders and that are also capable of handling the case of I.F.S. with {\em uncountably many} affine maps.
The proposed algorithms are based on the existence of a one-to-one transformation between the Jacobi matrix of $\mu$ and that of an auxiliary measure, $\sigma$, that encodes the parameters of affine, homogeneous I.F.S., that we will define momentarily.

Reversibility of such transformation also enables us to solve an inverse / approximation problem---that of computing the set of I.F.S. maps generating a given measure $\mu$. Original results \cite{steve} indicate
that I.F.S. quadratures, derived from the solution of the inverse problem, might offer
advantages with respect to conventional ones---yet, as was to be expected, numerical stability is an important issue also in this
inverse problem. Despite previous efforts \cite{ba0,ba2,steve,gio1,turch,ed} the first algorithm really achieving large orders in a stable way has been announced in \cite{nalgo1}: since it completes the theory developed herein and it can serve to prove experimentally the stability of the forward algorithms, it will be briefly discussed in the final section of this paper.

This paper is organized as follows: in the next section we review
the formalism of homogeneous iterated function systems and we generalize it to
allow for uncountably many maps. This leads to the definition of a convolution--like operator on measures, that is studied in Sect. \ref{secphi}, particularly with respect to its action on orthogonal polynomials and Jacobi matrices. This section is the core of the paper, since it contains almost all technical lemmas.
In Sect. \ref{secalgos} we use these lemmas to build numerical algorithms for the computation of the I.F.S. convolution and of the Jacobi matrix associated with an I.F.S., that are experimentally examined with respect to stability and performance in a number of significant cases. In Sect. \ref{secanal} these properties and the preceding theory allow us to derive rigorous results on, as well as numerical estimates of, significant analytical properties of I.F.S. measures such as their continuity and the capacity of their support.
In Sect. \ref{secspec} we describe a different technique for computing I.F.S. convolutions and Jacobi matrices,
based on discrete measures and on inverse Gaussian methods. To complete the paper, in Sect. \ref{secinve} we discuss the inverse algorithm, that yields the I.F.S. approximations of a target measure. The conclusions, Sect. \ref{secconclu}, briefly sum up the work and discuss further developments.

\section{Homogeneous Affine Iterated Function Systems}
\label{sechomo}

Systems of iterated functions \cite{hut,dem,ba2} are finite
collections of real maps $\phi_i : {\bf R} \rightarrow
{\bf R}$, $i = 1, \ldots, M$, for which there exists a
set ${\mathcal A}$, called the {\em attractor} of the I.F.S., that solves
the equation
 $
    {\mathcal A}=\bigcup_{i=1,\ldots ,M}\;\phi_i({\mathcal A}).
 $
Existence and uniqueness of ${\mathcal A}$ can be ensured under general circumstances. In this paper, we adopt a specific choice of the maps $\phi$ that has the advantage of leading to a
structured algebraic problem: that of {\em homogeneous, affine}
transformations:
\begin{equation}
\label{mappi}
    \phi_i (s) = \delta (s - \beta_i) + \beta_j  \;\;  i = 1, \ldots, M ,
\end{equation}
where $\delta$ is a real constant between zero and one, and
$\beta_{i}$ are real constants, that geometrically correspond to the fixed points of the maps. By
associating a positive weight, $\pi_{i} > 0$,
$\sum_{i} \pi_{i} = 1$, to each map one can define a measure $\mu$ supported on ${\mathcal A}$  via a procedure that we generalize in eq. (\ref{bala2}) below. The measure $\mu$ is specified by
the value of $\delta$ and of the pairs $(\beta_i,\pi_i)$, for $i =
1, \ldots, M$. This completes the definition of a ``classical'' homogeneous I.F.S.

Rather than restricting the cardinality of maps, $M$, to be finite (or countable, as done by Maulin and Urbansky \cite{urba}) we now allow the index set to be a continuum. This can be done in a variety of ways \cite{mendiv}. Our approach is to follow \cite{elton} and to observe that a finite,  homogeneous I.F.S. is fully described by the choice of $\delta$ and of the discrete measure
\begin{equation}
\label{sigma1}
   \sigma = \sum_{j=1}^M \pi_j D_{\beta_j},
\end{equation}
where $D_{x}$ is a unit mass atomic (Dirac) measure at the point $x$.
We now let any positive probability measure
$\sigma$ to be the distribution of affine constants: we only assume that the support of $\sigma$ is contained in a finite interval, which, without loss of generality, we may take to be $[-1,1]$.

\begin{definition}
\label{gafhom}
Let $\sigma$ be a positive Borel probability measure
on ${\bf R}$ whose support is contained in $[-1,1]$,
let $\delta$ be a real  number in $[0,1)$ and let $\bar{\delta} := 1 - \delta$.  Let the real number $\beta$ parameterize the I.F.S. maps $\phi(\beta,\cdot)$ as
$
 \phi(\beta,s) := \delta s + \bar{\delta} \beta.
$
The invariant I.F.S.
measure associated with the {\em affine homogeneous
I.F.S.} $(\delta,\sigma)$ is the unique probability measure $\mu$ that satisfies
\begin{equation}
\label{bala2}
    \int   f(s) \; d\mu(s)  =    \int d \sigma(\beta) \int d\mu(s) \;
    f (\phi(\beta,s)),
\end{equation}
for any continuous function $f$.
\end{definition}

\begin{remark}
If $\delta=0$, then  $\mu=\sigma$.  In fact, in this case
$\phi(\beta,s) = \beta$ for all $s$.
\end{remark}

Consistency of Def. \ref{gafhom}, {\em i.e.} existence and
uniqueness of $\mu$ are easily proven to hold:
\begin{proposition}
For any $\delta \in [0,1)$ eq. (\ref{bala2}) defines an invertible
transformation from the space ${\mathcal M}([-1,1])$ of probability
measures $\sigma$ on $[-1,1]$ to the subset of ${\mathcal M}([-1,1])$ composed of invariant measures $\mu$ of homogeneous I.F.S. with contraction ratio $\delta$. The support of $\mu$ contains the support of $\sigma$ and the convex hulls of these two sets coincide.
\label{lemexiste1}
\end{proposition}
{\em Proof.}
This result is already contained in \cite{elton}, although not entirely, because of a different parametrization of I.F.S. maps.
Let us start from the second statement. Recall that the attractor ${\mathcal A}$ of an I.F.S., which is the support of the invariant measure $\mu$, can also be characterized as the closure of the set of fixed points of the composition of an arbitrary number of maps $\phi(\beta,\cdot)$. It is then apparent that the support of $\sigma$ is a subset of that of $\mu$, and that the convex hull of the support of $\mu$ is the interval between the infimum and the supremum of the set of fixed points $\beta$, that is, the convex hull of the support of $\sigma$.

Next, since the supports of $\sigma$ and $\mu$ are contained in a finite interval, their moment problem is determined \cite{ach} and their infinite sets of moments uniquely identify the corresponding measures. The two sets of moments can be put in one-to-one relation: putting $f(s)=s^j$ in eq. (\ref{bala2}), for any $j$, one gets
a triangular set of relations by which it is
possible to derive all moments of $\mu$ from those of $\sigma$,
and viceversa. \qed

\begin{remark}
The proof above does {\em not} offer a numerically viable mean of computing the described correspondence between $\mu$ and $\sigma$, at least as far as Jacobi matrices are concerned, as noted in the Introduction.
\end{remark}

\section{A convolution-like operator and orthogonal polynomials}
\label{secphi}
We now want to define a mapping $\Phi$ in the space ${\mathcal M}([-1,1])$ of probability measures on $[-1,1]$,
according to which the invariant measure $\mu$ of an I.F.S. is the fixed point of such
transformation: $\mu = \Phi (\mu)$. Since $\Phi$ turns out to be contractive in a suitable metric, $\mu$ can be found as the
limit of the sequence of measures $\mu_n := \Phi^n \mu_0$, where
$\mu_0$ is any initial probability measure. To achieve this goal, we first describe a convolution-like operator of measures induced by the I.F.S. construction and derive its basic analytical properties (next subsection). We then describe the scaling relations of orthogonal polynomials with respect to affine transformations (second subsection). Thanks to these relations, the action of the operator $\Phi$ can be finally transferred on the Jacobi matrices of measures in ${\mathcal M}([-1,1])$ (third subsection).

\subsection{I.F.S. convolution of measures}
\label{ifsconvo}

We define a convolution like-operator $\Phi_\delta(\cdot;\cdot)$ acting on ${\mathcal M}^2([-1,1])$ as follows.

\begin{definition}
\label{defconv}
Let $\sigma, \eta$ two measures in ${\mathcal M}([-1,1])$ and let $\delta \in [0,1)$. The measure $\bar{\eta} \in {\mathcal M}([-1,1])$, called the IFS convolution of $\sigma$ and $\eta$,  $\bar{\eta} := \Phi_\delta(\sigma;\eta)$, is defined via the equation
\begin{equation}
\label{bala2x}
    \int   f(s) \; d \bar{\eta}(s)  =    \int d \sigma(\beta) \int d\eta (s) \;
        f (\phi(\beta,s)),
    \end{equation}
holding for any continuous function $f$.
\end{definition}
\begin{remark}
The symmetric r\^ole of $\sigma$ and $\eta$ in eq. (\ref{bala2x}) shows that they can be interchanged in the action of $\Phi_\delta(\cdot;\cdot)$, provided one also exchanges $\delta$ and $\bar{\delta}$:
 $    \Phi_\delta(\sigma;\eta) = \Phi_{\bar{\delta}}(\eta;\sigma)$.
\end{remark}

Suppose now that one keeps $\sigma$ fixed in the above construction. Then, $\bar{\eta} := \Phi_\delta(\sigma;\eta)$ is a function of $\eta$ alone: this defines the mapping
$\Phi_\delta(\sigma;\cdot)$ in ${\mathcal M}([-1,1])$. According to well established theory, see {\em e.g.} \cite{elton,mendiv}, one can easily prove
\begin{proposition}
The mapping $\Phi_\delta(\sigma;\cdot)$ defines a contraction in ${\mathcal M}([-1,1])\;$  equipped with the Hutchinson metric.
\label{propcontra}
\end{proposition}
Clearly, $\Phi_\delta(\sigma;\cdot)$ also
induces a map between the Jacobi matrices associated
with the related measures. We want now to derive the algebraic
relations required to translate this mapping into a numerically stable procedure.

\subsection{Scaling properties of orthogonal polynomials}
\label{subsescaling}

Fundamental for the algebraic structure of the algorithms that we are about to develop is the
study of the scaling properties of polynomials
with respect to the maps $\phi$. We start with the
almost trivial, but fundamental

\begin{lemma}
 \label{lem3p1}
Let $u_n(s)$ be any polynomial of degree $n$. Let
$\{p_n(\eta;s)\}_{n \in \bf N}$  and $\{p_n(\sigma;\beta)\}_{n \in
\bf N}$ be the families of orthogonal polynomials associated with the
measures $\eta$ and $\sigma$ in the variables $s$ and $\beta$,
respectively. Then, there exists constants $\Omega^n_{k,r}$, with $k,r \geq 0$, such that
 \beq
    u_n(\phi(\beta,s)) =
     \sum_{0 \leq k+r \leq n }  \Omega^n_{k,r}
     p_k(\eta; s)
 p_r(\sigma; \beta).
 \label{gamma1}
 \end{equation}
 \label{lemgamma}
\end{lemma}
{\em Proof.}
The polynomial $u_n(\phi(\beta,s)) = u_n(\delta s + \bar{\delta}
\beta)$ can be written in the form $\sum_{j=0,\ldots,n} c_{j} s^j
\beta^{n-j}$. Expanding the monomials $s^j$ and $\beta^{n-j}$ in
orthogonal polynomials of $\eta$ and $\sigma$ respectively proves eq.
(\ref{gamma1}). \qed

In this paper, we shall use Lemma \ref{lemgamma} specifically for
$u_n(s)=p_n(\bar{\eta};s)$, the n-th orthogonal polynomial of the measure $\bar{\eta} =  \Phi_\delta(\sigma;\eta)$ in  eq. (\ref{bala2x}). $\Omega^n$, with $n=0,1,\ldots$, 
will designate uniquely the related coefficients and will
be called the {\em scaling  matrices}. For simplicity of notation,
we shall not indicate explicitly the $\bar{\eta},\eta,\sigma$ dependence of the matrix $\Omega^n$. Also,
we shall let the indices $k$ and $r$ run freely, while assuming that
$\Omega^n_{k,r} = 0$ unless $k,r \geq 0$ and $0 \leq k+r \leq n$, so that effectively $\Omega^n$ is a triangular matrix.
The notation $q_n$ will indicate a polynomial of degree $n$ of which
no further specification is necessary.

\begin{lemma}
For any $n \in {\bf N}$, the scaling matrix $\Omega^n$ satisfies the
normalization condition
 \beq
   \sum_{k,r}
  (\Omega^n_{k,r})^2 = 1.
 \label{gam34}
 \end{equation}
 \label{lemgnorm}
\end{lemma}
{\em Proof.}
Since $\int d\bar{\eta}(s) p_{n}(\bar{\eta};s)^2 = 1$, eq. (\ref{bala2x}) implies
that
 \[
  1 =
  \dbleint d \sigma(\beta)   d\eta(s)
   p_{n}(\bar{\eta};\phi(\beta,s))^2 =
  \]
 \[
   =
  \dbleint d \sigma(\beta)   d\eta(s)
  \! \! \! \! \!
\sum_{k,r,k',r'} \! \! \! \!
  \Omega^n_{k,r} \Omega^n_{k',r'}
 p _{k}(\eta; s)
 p_{r}(\sigma; \beta)
  p _{k'}(\eta; s)
 p_{r'}(\sigma; \beta).
\]
Performing the integrations with the aid of orthogonality leads to
eq. (\ref{gam34}). \qed

\begin{lemma}
For any $n \in {\bf N}$, the scaling  matrix extremal entries $\Omega^n_{n,0}$ and
$\Omega^n_{0,n}$ are given by
\beq
 \Omega^n_{n,0} = \frac{h_n(\bar{\eta})}{h_n(\eta)} \delta^n, \mbox{ and } \;
    \Omega^n_{0,n} = \frac{h_n(\bar{\eta})}{h_n(\sigma)} \bar{\delta}^n,
    \label{ga31e}
      \end{equation}
where $h_n(\cdot)$ denote the coefficients of the highest power in
the orthonormal polynomial of degree $n$: $p_n(\cdot;s) = h_n(\cdot) s^n + q_{n-1}(s)$.
 \label{lg1e}
 \end{lemma}
{\em Proof.}
Clearly, because of eq. (\ref{gamma1}), and because
$p_0(\sigma;\beta) = 1$,
on
the one hand
\[
  p_n(\bar{\eta};\delta s + \bar{\delta} \beta) = \Omega^n_{n,0} p_n(\eta;s) +
    q_{n-1}(s) = \Omega^n_{n,0} h_n(\eta) s^n +
      q'_{n-1}(s),
        \]
where $q_{n-1}(s)$ and $q'_{n-1}(s)$ are polynomials of degree
$n-1$ in the variable $s$ and the $\beta$ dependence has been implied. On the other hand,
\[
 p_n(\bar{\eta};\delta s + \bar{\delta} \beta) =  h_n(\bar{\eta}) \delta^n s^n +
  q''_{n-1}(s),
   \]
with $q''_{n-1}(s)$ another polynomials of degree $n-1$. This proves the first  eq. (\ref{ga31e}). Because of the symmetrical r\^ole of $\sigma$ and $\eta$, this also proves the second. \qed

\begin{lemma}
\beq
  \phi(\beta,s) p_n(\bar{\eta};\phi(\beta,s)) =
\sum_{k,r }  \Omega^n_{k,r}
  [ \delta P_k(\eta;s)p_r(\sigma; \beta)
  + \bar{\delta} P_r(\sigma;\beta)p _k(\eta; s) ],
 \nuq{funda1e}
where we have put
  $
P_k(\cdot;s) :=
       b_{k+1}(\cdot) p _{k+1}(\cdot; s) +
a_{k}(\cdot) p _{k}(\cdot; s)  + b_{k}(\cdot) p _{k-1}(\cdot; s)
$
and $\cdot$ can be either $\eta$ or $\sigma$. \label{lemgamm2}
\end{lemma}
{\em Proof.}
Because of eq. (\ref{gamma1})
$$
 \phi(\beta,s) p_n(\bar{\eta};\phi(\beta,s)) =
  (\delta s +
\bar{\delta} \beta)
  \sum_{k,r }  \Omega^n_{k,r}
     p _k(\eta; s)
 p_r(\sigma; \beta)
 $$
and the products $s p_k(\eta;s)$ and $\beta p_r(\sigma;\beta)$ can be
dealt with using the recurrence relation for orthogonal polynomials,
eq (\ref{nor2}). \qed

\subsection{Relations between Jacobi and scaling  matrices}
We can now show that the scaling properties embodied in the
matrices $\Omega^n$ imply algebraic relations between the Jacobi matrices
$J_\eta$, $J_{\bar{\eta}}$ and $J_\sigma$. In other words, we are capable of translating on these latter the action of the convolution operator: $J_{\bar{\eta}} = \Phi_\delta(J_\sigma;J_\eta)$, where, with a slight abuse of notation, we also denote by $\Phi_\delta(\cdot;\cdot)$ the induced operator in the space of pairs of Jacobi matrices. 

\begin{lemma}
For any $n \in {\bf N}$ and $\delta \in [0,1)$, the Jacobi matrix
entry $a_n(\bar{\eta})$ can be written as a linear combination of
the Jacobi matrix entries
$
  \{a_j(\eta),b_j(\eta) \}_{j=0}^n
   $
   and
   $
  \{a_j(\sigma),b_j(\sigma) \}_{j=0}^n,
   $
   with coefficients derived from the scaling  matrix
$\Omega^n$:
\begin{equation}
 a_n (\bar{\eta})  =
 \sum_{k,r} \Omega^n_{k,r}
   \left(
   \delta [ a_k(\eta)  \Omega^n_{k,r}  + 2 b_k(\eta)
  \Omega^n_{k-1,r}] +
  \bar{\delta} [ a_r(\sigma)  \Omega^n_{k,r}  + 2 b_r(\sigma)
  \Omega^n_{k,r-1}] \right).
\label{ad2d}
\end{equation}
The coefficients of the highest index terms $a_n(\eta)$ and
$a_n(\sigma)$ at r.h.s. are ${\delta} (\Omega^n_{n,0})^2$ and $\bar{\delta}
(\Omega^n_{0,n})^2$, respectively.
  \label{lemgamm3}
\end{lemma}
{\em Proof.}
We use eq. (\ref{nor2}) and orthogonality to write
 \beq
 a_n (\bar{\eta}) = \int d\bar{\eta}(s) s p^2_n(\bar{\eta};s)
 \nuq{ad1}
and we then expand according to the relation (\ref{bala2x})
and to eqs. (\ref{gamma1}), (\ref{funda1e}):
 \beq
  a_n (\bar{\eta}) =  \int d \sigma(\beta) \int d\eta(s) \phi(\beta,s)
  [p_n(\bar{\eta};\phi(\beta,s))]^2 =
\label{ad2a}
\end{equation}
 $$ = \dbleint d \sigma(\beta)   d\eta(s)
  \! \! \! \! \!
\sum_{k,r,k',r'} \! \! \! \!
  \Omega^n_{k,r} \Omega^n_{k',r'}
   p _{k'}(\eta; s)
 p_{r'}(\sigma; \beta)
   [ \delta P_k(\eta;s) p_{r}(\sigma; \beta)+ \bar{\delta} P_r(\sigma;\beta)  p _{k}(\eta;
   s)].
$$
Integrations can be computed explicitly using orthogonality of the
polynomials: after some manipulations, one gets the linear relation
eq. (\ref{ad2d}).
Notice that the highest index with non-zero coefficient for both $k$
and $r$ is $n$.
Direct inspection shows that the coefficient of $a_n(\eta)$ at
r.h.s. is ${\delta} (\Omega^n_{n,0})^2$ and that of $a_n(\sigma)$
is $\bar{\delta} (\Omega^n_{0,n})^2$. \qed


\begin{lemma}
\label{lgammarec}
For any $n \in {\bf N}$, all the entries of the matrix
$\tilde{\Omega}^{n+1} := b_{n+1}(\bar{\eta}) \Omega^{n+1}$ can be
computed linearly in terms of the entries of $\Omega^{n}$ and
$\Omega^{n-1}$, with coefficients determined by the Jacobi matrix
entries $\{a_j(\eta)\}_{j=0}^n$, $\{a_{j}(\bar{\eta})\}_{j=0}^n$,
$\{a_{j}(\sigma)\}_{j=0}^n$, $\{b_{j}(\bar{\eta})\}_{j=0}^n$,
 $\{b_j(\eta)\}_{j=1}^{n+1}$,
and $\{b_j(\sigma)\}_{j=1}^{n+1}$. In fact, the following relation holds:

$$
    \tilde{\Omega}^{n+1}_{j,l} := b_{n+1}(\bar{\eta})
\Omega^{n+1}_{j,l} =
 \delta \left[ a_j(\eta) \Omega^{n}_{j,l} + b_{j+1}(\eta) \Omega^{n}_{j+1,l}
  + b_{j}(\eta) \Omega^{n}_{j-1,l} \right]
 +
  $$
 \beq
 +  \bar{\delta} \left[ a_l(\sigma) \Omega^{n}_{j,l} + b_{l+1}(\sigma)
\Omega^{n}_{j,l+1}
  + b_{l}(\sigma) \Omega^{n}_{j,l-1} \right]
   - a_n(\bar{\eta}) \Omega^{n}_{j,l}
   - b_n(\bar{\eta}) \Omega^{n-1}_{j,l}.
\label{gammarec}
\end{equation}

The entries $b_{n+1}(\sigma)$, $b_{n+1}(\eta)$, only affect the
computation of the extremal values
 \beq
  \tilde{\Omega}^{n+1}_{n+1,0} :=
b_{n+1}(\bar{\eta}) \Omega^{n+1}_{n+1,0} = \delta \Omega^{n}_{0,n}
b_{n+1}(\eta),
 \nuq{lasterm2}
 \beq
  \tilde{\Omega}^{n+1}_{0,n+1} :=
b_{n+1}(\bar{\eta}) \Omega^{n+1}_{0,n+1} = \bar{\delta} \Omega^{n}_{0,n}
b_{n+1}(\sigma).
 \nuq{lasterm}
\end{lemma}
{\em Proof.}
We start from the relation
  \beq
 b_{n+1}(\bar{\eta}) p_{n+1} (\bar{\eta};s) =
(s-a_n(\bar{\eta})) p_n(\bar{\eta};s) - b_n(\bar{\eta})
p_{n-1}(\bar{\eta};s),
 \nuq{use2}
 and map $s$ in
$\phi(\beta,s) = \delta s + \bar{\delta} \beta$. Define $Q_{n+1}(s)
:= b_{n+1}(\bar{\eta}) p_{n+1} (\bar{\eta};\phi(\beta,s))$:
 \beq
 Q_{n+1}(s) =
(\phi(\beta,s) -a_n(\bar{\eta}) ) p_n(\bar{\eta};\phi(\beta,s))
  - b_n(\bar{\eta}) p_{n-1}(\bar{\eta};\phi(\beta,s)).
 \nuq{pigra1}
Use now eqs. (\ref{gamma1}) and (\ref{funda1e}) to express the terms
at r.h.s. via the matrices $\Omega^n$, and $\Omega^{n-1}$.
Then, we multiply both sides of eq. (\ref{pigra1}) by $p_{j}(\eta;s)
p_{l}(\sigma;\beta)$ and integrate w.r.t. $d \sigma(\beta)
d\eta(s)$, to get eq. (\ref{gammarec}).
The r.h.s. of this equation is a linear
combinations of the matrix elements of $\Omega^{n}$ and
$\Omega^{n-1}$. The coefficients are given by the matrix entries of
$J_{{\eta}}$, $J_{\bar{\eta}}$ and $J_\sigma$. Direct inspection,
using the triangular nature of the $\Omega^n$'s, reveals that the
terms of highest index appearing at r.h.s. are $a_n(\eta)$,
$b_{n+1}(\eta)$ and $a_n(\sigma)$, $b_{n+1}(\sigma)$, and it also reveals that the
last part of the thesis and eqs. (\ref{lasterm2}), (\ref{lasterm}) hold. \qed

\begin{lemma}
For any $n \in {\bf N}$, the Jacobi matrix entries $b_{n+1}(\bar{\eta})$,
$b_{n+1}({\eta})$ and $b_{n+1}(\sigma)$ are related to the
scaling  matrix $\tilde{\Omega}^{n+1}$ via
\beq
 \label{gam353}
 b^2_{n+1}(\bar{\eta})
   = {\delta}^2
    b^2_{n+1}(\eta)
  (\Omega^{n}_{n,0})^2
   +  \bar{\delta}^2 b^2_{n+1}(\sigma)
  (\Omega^{n}_{0,n})^2 +
   \sum_{k,r} ^\prime
  (\tilde{\Omega}^{n+1}_{k,r})^2,
 \end{equation}
where the primed summation runs over all pairs of indices $(k,r)$ that are different from $(n+1,0)$ and $(0,n+1)$.
 \label{lemrnp2}
\end{lemma}
{\em Proof.}
Because of eq. (\ref{gam34}),

 \beq
 \label{gam352}
 b^2_{n+1}(\bar{\eta}) =
 \sum_{k,r}
  (\tilde{\Omega}^{n+1}_{k,r})^2 =
(\tilde{\Omega}^{n+1}_{0,n+1})^2 +
  (\tilde{\Omega}^{n+1}_{n+1,0})^2 +
   \sum_{k,r}'
  (\tilde{\Omega}^{n+1}_{k,r})^2.
 \end{equation}
We can now use the explicit formulae (\ref{lasterm2}) and
(\ref{lasterm}). \qed

\section{Algorithms for the construction of the I.F.S. Jacobi Matrix}
\label{secalgos}
Having derived the necessary technical lemmas in the previous section, we
are now in a position to chain them into two algorithms for the construction of the Jacobi matrix $J_\mu$ of an affine, homogeneous I.F.S. described by the contraction ratio $\delta$ and by the affine constants distribution $\sigma$. We first develop a {\em fixed--point, forward algorithm} similar to that already exploited in \cite{mobius}, based on a technique for computing the I.F.S. convolution of two Jacobi matrices.
Next, the solution of the fixed point equation can be obtained by a bootstrap technique like that of ref. \cite{cap}. This yields the {\em closure algorithm}.
We will denote by $J^{(\bar{n})}$ the finite truncation of any Jacobi matrix $J$ to size ${\bar{n}}$.

\subsection{I.F.S. convolution of Jacobi matrices}

\begin{theorem} The Jacobi matrix $J_{\bar{\eta}}$ of the I.F.S. convolution measure
$\bar{\eta} := \Phi_\delta(\sigma;\eta)$ can be computed recursively from the
Jacobi matrices $J_\sigma$ and $J_\eta$. Computation of a finite
truncation $J^{(\bar{n})}_{\bar{\eta}}$  of size ${\bar{n}}$ of $J_{\bar{\eta}}$ only requires knowledge of the truncated matrices
$J_\sigma^{(\bar{n})}$ and $J^{(\bar{n})}_{{\eta}}$. \label{theodirect2}
\end{theorem}
{\em Proof.}
The algorithm is structured in the following sequence of steps
\begin{itemize}
 \item[] {{\bf Algorithm 1.\\ Computing the Jacobi matrix of an I.F.S. convolution.} \\
 {\bf Input}: the (truncated) Jacobi matrices $J_\sigma^{(\bar{n})}$  and $J^{(\bar{n})}_{{\eta}}$, the contraction factor $\delta$, the truncation size $\bar{n}$.
 \\ {\bf Output}: the (truncated) Jacobi matrix of  $\bar{\eta} := \Phi_\delta(\sigma;\eta)$.}
\item[0:] Initialization: $n=0$. One has $\Omega^0_{0,0}= 1$, since
$p_0(\bar{\eta};s)=p_0({\eta};s)=p_0(\sigma;\beta)=1$ and
$b_0(\bar{\eta})=b_0({\eta})=b_0(\sigma)=0$.
\item[1:] Induction hypothesis:
$\{\Omega^j,j=0,\ldots,n\}$, $\{a_j(\bar{\eta}), j=0,\ldots,n-1 \}$
and $\{b_j(\bar{\eta}), j=0,\ldots,n \}$ are known.
\item [2:] Computation of $a_n(\bar{\eta})$: eq. (\ref{ad2d}) in Lemma \ref{lemgamm3}.
\item [3:] Computation of the matrix $\tilde{\Omega}^{n+1}$: Lemma
\ref{lgammarec}.
\item [4:] Computation of $b_{n+1}(\bar{\eta})$: eq. (\ref{gam353}) in Lemma \ref{lemrnp2}
\item [5:] Computation of $\Omega^{n+1}$: divide
$\tilde{\Omega}^{n+1}$ by $b_{n+1}(\bar{\eta})$.
\item [6:] Augment $n$ to $n+1$ and loop back to step 1 if $n+1$ is less than the desired truncation $\bar{n}$, otherwise stop.
\end{itemize}
Notice that $b_{n+1}(\bar{\eta})$ in steps [4] and [5] is never zero, if $\delta>0$, or if the cardinality of the support of either $\sigma$ or $\eta$ is larger than $n+1$.
\qed

\subsection{Fixed--point Forward Algorithm}
The previous algorithm can serve a double purpose: on the one hand,
to investigate the convolution--like measure
$\bar{\eta}$ as a function of the factors $\sigma$ and $\eta$. On
the other hand, to obtain the Jacobi matrix of the invariant
$(\delta,\sigma)$--IFS measure in an iterative fashion, as in the following:

\begin{itemize}
 \item[] {{\bf Algorithm 1-Fix.    Computing the $  (\delta,\sigma)$--IFS Jacobi matrix.} \\
{\bf Input}: the (truncated) Jacobi matrix $J_\sigma^{(\bar{n})}$ of the distribution of affine constants $\sigma$, the contraction factor $\delta$, the truncation size $\bar{n}$ and a convergence threshold.
 \\ {\bf Output}: the (truncated) Jacobi matrix $J_\mu^{(\bar{n})}$ of the I.F.S. measure $\mu$.}
\item[0:] Initialization: Let $\mu_0$ be a positive probability
measure in ${\cal M}([-1,1])$ and let $J_0^{(\bar{n})}$ be its (truncated) Jacobi matrix.
\item[1:] For $m=1$ until convergence \\
Use algorithm 1 to compute $J^{(\bar{n})}_m$, the truncated Jacobi matrix of
$\mu_m:=\Phi_\delta(\sigma;\mu_{m-1})$.
\end{itemize}

\begin{remark}
Convergence of the previous algorithm is assured by the contractive nature of the transformation $\Phi_\delta(\sigma;\cdot)$, Prop. \ref{propcontra}.
Numerical convergence, on the other hand, can be gauged {\em e.g.} according to the Frobenius norm of the
difference $J^{(\bar{n})}_m - J_{m-1}^{(\bar{n})}$. Observe that, because of theorem \ref{theodirect2}, this difference is exact (except for numerical errors) at any finite truncation $\bar{n}$, that is to say, enlarging $\bar{n}$ does not change the values of the already computed smaller truncations of the Jacobi matrix.
\end{remark}

\subsection{Closure Algorithm}

We now show that the Jacobi matrix of the invariant measure
$\mu$ of a $(\delta,\sigma)$--IFS can be obtained directly and in a finite number of steps---that is, the sequence of operations of Algorithm 1-Fix can be exactly closed. In fact, we
can set $\mu=\eta=\bar{\eta}$ throughout the formulae of
section \ref{secphi}. The only notable difference occurs in the following lemma:

\begin{lemma}
For any $n$, the Jacobi matrix entries $b_{n+1}(\mu)$ and
$b_{n+1}(\sigma)$ are related to the scaling  matrix
$\tilde{\Omega}^{n+1}$ via
 \beq
 \label{gam356}
 (1 - \delta^{2(n+1)}) b^2_{n+1}(\mu)
   =
   b^2_{n+1}(\sigma) \bar{\delta}^2
  (\Omega^{n}_{0,n})^2 +
   \sum'_{k,r}
  (\tilde{\Omega}^{n+1}_{k,r})^2,
 \end{equation}
 where the primed summation runs over pairs of indices $(k,r)$ not
equal to $(n+1,0)$ and $(0,n+1)$. Therefore, the Jacobi matrix entry
$b_{n+1}(\mu)$ can be always computed from the other quantities in
the above equation, while $b_{n+1}(\sigma)$ is only defined from the other terms when the
difference between the l.h.s. and the second term at r.h.s. is
positive.
 \label{lemrnp2x}
\end{lemma}
{\em Proof.} This is a consequence of Lemma \ref{lemrnp2} and can be proven similarly. \qed

\begin{remark}
Observe the breaking of the symmetry in the r\^ole of the measures
$\mu$ and $\sigma$ in the previous lemma. This is of paramount
importance for the inverse problem discussed in Sect. \ref{secinve}.
\end{remark}

\begin{theorem} The Jacobi matrix $J_\mu$ of a homogeneous I.F.S.
with contraction ratio $\delta$ and affine constants distribution
$\sigma$ can be computed recursively from the Jacobi matrix
$J_\sigma$. Computation of a finite truncation $J_\mu^{(\bar{n})}$ only
requires knowledge of $J_\sigma^{(\bar{n})}$. \label{theodirect}
\end{theorem}
{\em Proof.}
The algorithm is effected in the following sequence of steps
\begin{itemize}
 \item[] {{\bf Algorithm 2.  Computing the $  (\delta,\sigma)$--IFS Jacobi matrix.}  \\
{\bf Input}: the (truncated) Jacobi matrix $J_\sigma^{(\bar{n})}$, the contraction factor $\delta$, the truncation size $\bar{n}$.
\\ {\bf Output}: the (truncated) Jacobi matrix  $J_\mu^{(\bar{n})}$ of the I.F.S. measure $\mu$}
\item[0:] Initialization: $n=0$. One has $\Omega^0_{0,0}= 1$, since
$p_0(\mu;s)=p_0(\sigma;\beta)=1$ and $b_0(\mu)=b_0(\sigma)=0$.
\item[1:] Induction hypothesis:
$\{\Omega^j,j=0,\ldots,n\}$, $\{a_j(\mu), j=0,\ldots,n-1 \}$ and
$\{b_j(\mu), j=0,\ldots,n \}$ are known.
\item [2:] Computation of $a_n(\mu)$: Lemma \ref{lemgamm3}. Observe
that setting $\eta=\bar{\eta}=\mu$ in eq. (\ref{ad2d}) the term
$a_n(\mu)$ appears at r.h.s. with a coefficient ${\delta}
(\Omega^n_{n,0})^2$ that is strictly less than one.
\item [3:] Computation of the matrix $\tilde{\Omega}^{n+1}$: Lemma \ref{lgammarec}.
\item [4:] Computation of $b_{n+1}(\mu)$: Lemma \ref{lemrnp2x}
\item [5:] Computation of $\Omega^{n+1}$: divide
$\tilde{\Omega}^{n+1}$ by $b_{n+1}(\mu)$.
\item [6:] Augment $n$ to $n+1$ and loop back to step 1 if $n$ is less than the desired truncation $\bar{n}$, otherwise stop.
\end{itemize}
Notice again that $b_{n+1}(\mu)$ in steps [4] and [5] is never zero, if $\delta>0$, or if the cardinality of the support of $\sigma$ is larger than $n+1$.
\qed

\begin{remark}
The algorithm can be carried out in ${\mathcal O} (\bar{n}^3)$ operations.
In addition, it can also be re-structured in order to compute in
place of the $b_n(\mu)$ the squares $b_n^2(\mu)$, that enter the
recursion relation of {\em monic} orthogonal polynomials ({\em i.e.}
those normalized as having unit coefficient in the highest power).
This avoids the square root in step 4 and therefore leads to an
algorithm that can be carried out exactly in rational arithmetics, when $\delta$ and $J_\sigma$ are such, or symbolically, for instance to obtain the recursion coefficients as functions of the contraction ratio $\delta$.
\end{remark}

\begin{remark}
The storage requirement of the algorithm is ${\mathcal O} (\bar{n}^2)$.
\end{remark}

\begin{remark}
In the classical I.F.S. case with $M$ maps, when $\sigma$ is a
finite sum of $M$ atomic measures, eq. (\ref{sigma1}), the algorithm
above can be used without any modification. In this case
the Jacobi matrix $J_\sigma$ is finite and one can set $b_j(\sigma)
= 0$ for $j \geq M$. This entails, via eq. (\ref{gammarec}), that $
         \Omega^{n}_{r,k} = 0
$ for all $k > M$. The algorithm then runs in ${\mathcal O} (M \bar{n}^2)$
operations and requires a storage of size ${\mathcal O} (M \bar{n})$,
exactly as the Stieltjes algorithm in \cite{cap}. Indeed, the
difference between the latter and the present algorithm is two--fold. Firstly, notice that in \cite{cap} 
the polynomial $p_n(\phi(\beta_j,s))$, for any $j=1,\ldots,M$, is developed on the basis
$\{p_k(\mu;s)\}_{k=0}^n$, resulting in $M$ coefficient vectors.
Here, we exploit the algebraic properties of these coefficients in
order to write them on the basis of the orthogonal polynomials of
$\sigma$.
Secondly, the input data of the algorithms are different: in \cite{cap} these are the parameters of a finite set of maps---here, they are the common contraction coefficient $\delta$ and the Jacobi matrix of the measure $\sigma$.
\end{remark}

\subsection{Numerical examples}
\label{sec-numex}
\begin{figure}[tb]
\includegraphics[width=6cm,height=12cm,angle=270]{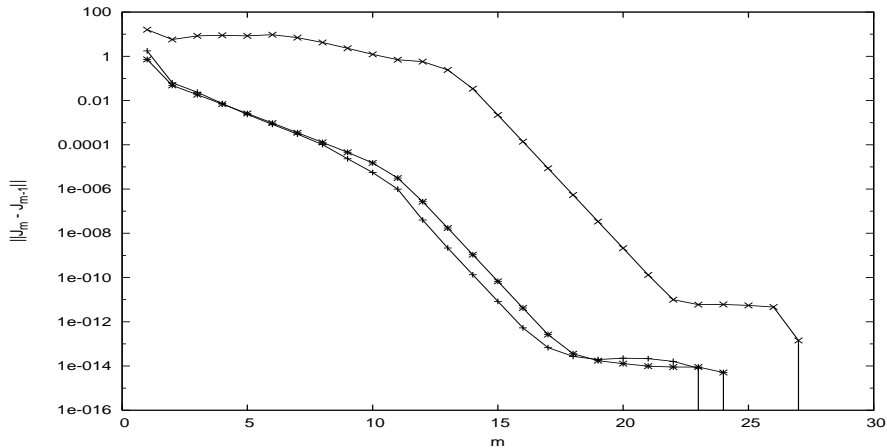}
\caption{Frobenius distance between the truncated Jacobi matrices of
$\mu_m$ and $\mu_{m-1}$, versus $m$, for $\bar{n}=4096$. The three measures $\sigma$ described in Examples \ref{exa-mupp} to \ref{exa-muac} are considered. The highest curve is for the atomic $\sigma_{pp}$ (Ex. \ref{exa-mupp}, crosses), while the two lowest sets of data are for $\sigma_{sc}$ (Ex. \ref{exa-musc}, pluses) and $\sigma_{ac}$ (Ex. \ref{exa-muac}, asteriscs).}\label{diffe4.fig}
\end{figure}

Let us start by examining the fixed--point, forward algorithm 1-Fix, choosing
$\mu_0$ as the uniform, normalized Lebesgue measure on $[-1,1]$.
This choice is dictated by the simplicity of its Jacobi matrix, but other choices would work equally well. Indeed, one might even start from a suitable Jacobi matrix whose associated measure is not even known. On the other hand, three kind of possible choices for $\sigma$ are:
\begin{itemize}
 \item  $\sigma_{pp}$:  a point measure with a finite number of atoms,
 \item $\sigma_{sc}$: a singular continuous measure,
 \item $\sigma_{ac}$: an absolutely continuous measure.
\end{itemize}
We investigate these different choices because of their r\^ole in mathematical physics, as described in the introduction and because of the open problems on their continuity properties, discussed in Sect. \ref{secanal}.
\begin{example}
Let the atomic measure $\sigma_{pp}
= \frac{1}{2} (D_{-1}+D_{1})$ be made of two atoms at positions
$-1$ and $1$, with equal weights. Fixing $\delta=\frac{3}{10}$
therefore implies that $\mu = \lim_{m \rightarrow \infty} \mu_m$ is the
invariant measure of a two-maps homogeneous I.F.S. with
$\delta=\frac{3}{10}$, equal weights $\pi_1=\pi_2=\frac{1}{2}$ and
$\beta_1=-\beta_2=1$. The convex hull of the support of this fractal
measure is also $[-1,1]$.
\label{exa-mupp}
\end{example}
\begin{example}
As to the second choice, $\sigma_{sc}$ can be taken as precisely
the two-maps I.F.S. measure of example \ref{exa-mupp}. The generated measure $\mu$ (with a value of the contraction ratio $\delta$ independent of that of Ex. \ref{exa-mupp}) is now the invariant measure of an I.F.S. with uncountably many maps, whose
fixed points are located on a Cantor set. For the case of Fig. \ref{diffe4.fig} we have chosen $\delta=1/4$.
\label{exa-musc}
\end{example}
\begin{example}
Thirdly, $\sigma_{ac}$ can be taken as the uniform Lebesgue measure
on $[-1,1]$. The generated measure $\mu$ (with $\delta=1/4$) is now an absolutely continuous measure supported on  $[-1,1]$, whose analytical properties have been studied in \cite{nalgo2}, where a graph of its density is also displayed.
\label{exa-muac}
\end{example}
We exhibit the numerical convergence properties of algorithm 1-Fix, when applied to Examples \ref{exa-mupp} to \ref{exa-muac}. In Fig. \ref{diffe4.fig} we plot, in log-linear scale, the Frobenius norms of the differences between the Jacobi matrices of $\mu_m$ and $\mu_{m-1}$, versus the iteration number $m$, computed at a fixed finite truncation. After an initial transient, that lasts longer for the point measure than in the other cases, we observe the exponential convergence typical of fixed point techniques. Obviously, the Jacobi matrix of the three measure just considered can be also computed via the closure algorithm 2, with an obvious saving of computer time. In certain investigations, though, the recursive algorithm may be needed, when one wants to study the properties of the Jacobi matrices of a sequence of measures converging to the limit measure $\mu$.

Let us now consider Algorithm 2.
In the case when the measure $\sigma$ is composed of a finite number of atoms, both storage and cpu time can be greatly reduced, as noted above. This fact permits far-reaching numerical experimentations. We first consider the uniform Lebesgue measure on $[-1,1]$, that can be generated by choosing  $\delta=1/2$ in example \ref{exa-mupp} while leaving $\sigma$ unchanged. Being the Jacobi matrix explicitly known, we have tested the error propagation, finding errors less that $6 \cdot 10^{-16}$ 
for $n$ as large as $250,000$, being the machine epsilon of the order of $2.22 \; 10^{-16}$.

Next, we consider the refinable functionals introduced in Ref. \cite{dirkvil}. As a matter of facts, their invariant measures are supported on a different interval than $[-1,1]$, but our algorithms work equally well without change.

\begin{example}
In I.F.S. language, the first numerical example discussed in Sect. 4 of Ref. \cite{dirkvil} consists of four maps, with contraction factor $\delta=1/2$ and with fixed points $\beta_j= j$, for $j=0,\ldots,3$ and weights $1/8, 3/8, 3/8, 1/8$, respectively.
\label{exa-dirk1}
\end{example}

\begin{figure}[tb]
\centerline{\includegraphics[width=6cm,height=12cm,angle=270]{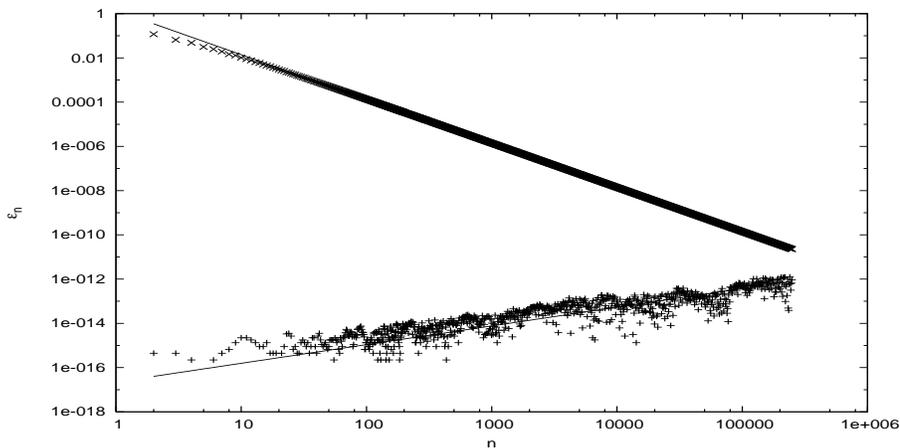}}
\caption{Absolute differences $\epsilon_n = |a_n-3/2|$ (lower data points, pluses) and $\epsilon_n = |b_n-3/4|$ (higher data points, crosses) versus $n$ for the refinable functional defined in Example \ref{exa-dirk1}. The fitting power-laws $\epsilon_n \sim n^\gamma$ have exponent $\gamma=.85$ and $\gamma=-2$, respectively }\label{dirk1a.fig}
\end{figure}

We have computed the Jacobi matrix of Example \ref{exa-dirk1} with algorithm 2 up to $n=250,000$ without encountering any numerical instability. Because of symmetry, one can theoretically assess that $a_n=3/2$, a value that is not automatically reproduced by the algorithm so that it can be used to gauge numerical error propagation. In figure \ref{dirk1a.fig} we observe a mild, slower than linear error growth for the diagonal coefficients, that is orders of magnitude lower than the difference between $b_n$ and the asymptotic limit $b_\infty=3/4$, also reported in the figure. The observed power-law decay of such  difference, with exponent $\gamma=-2$, is therefore to be deemed reliable and suggestive of the absolute continuity of the orthogonality measure $\mu$, according to the theory presented in Sect. \ref{secanal}.

\begin{figure}[tb]
\centerline{\includegraphics[width=6cm,height=12cm,angle=270]{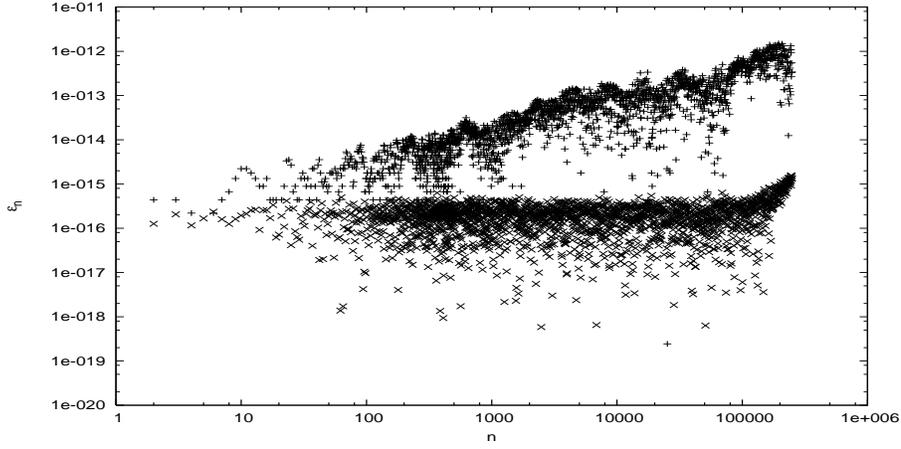}}
\caption{Absolute differences $\epsilon_n = |a_n-\tilde{a}_n|$ (crosses) and $\epsilon_n = |b_n-\tilde{b}|$ (pluses)  versus $n$ between the $a_n,b_n$  provided by algorithm 2 and the $\tilde{a}_n,\tilde{b}_n$ from the algorithm of Ref. \cite{cap} when applied to the same case of Ex. \ref{exa-dirk1}, Fig. \ref{dirk1a.fig}.} \label{dirk1b.fig}
\end{figure}

To provide a more stringent verification, we have also computed, in the same case as above, the difference between the results of Algorithm 2 and those of the Stieltjes technique of Ref. \cite{cap}, this latter run in quadruple precision (machine epsilon of the order of $10^{-34}$). The slower than linear growth of the absolute differences in the diagonal and the outdiagonal components of the Jacobi matrix, reported in Fig. \ref{dirk1b.fig}, suggests the numerical stability of both techniques.

\begin{figure}[tb]
\centerline{\includegraphics[width=6cm,height=12cm,angle=270]{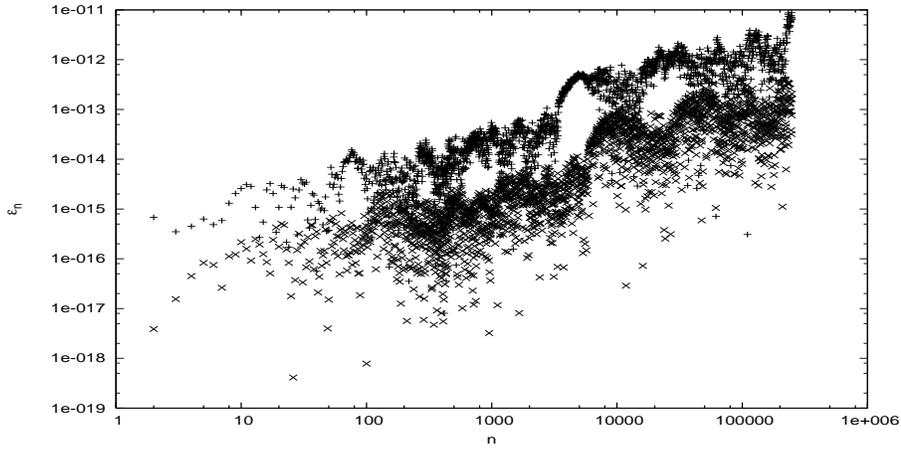}}
\caption{Absolute differences $\epsilon_n = |a_n-\tilde{a}_n|$ (crosses) and $\epsilon_n = |b_n-\tilde{b}|$ (pluses)  versus $n$ between the $a_n,b_n$  provided by algorithm 2 and the $\tilde{a}_n,\tilde{b}_n$ from the algorithm of Ref. \cite{cap} when applied to Ex. \ref{exa-dirk2}.} \label{dirk2.fig}
\end{figure}

Finally, we have run the same test on the second example proposed by Laurie \cite{dirkvil}, who observed a linear growth of the error in his technique, discovering a comparable behavior, on the larger range here displayed: Ex. \ref{exa-dirk2} and Fig. \ref{dirk2.fig}. Using the theory of the next section, absolute continuity of the associated measure can also be conjectured from the computed $J_\mu$.
\begin{example}
In example \ref{exa-dirk1}, use the modified weights $1/8, 3/8, 3/8, 1/8$.
\label{exa-dirk2}
\end{example}

We propose numerical experiments on the case of infinite Jacobi matrices $J_\sigma$ later on, in Sect. \ref{secinve}. We now pause for a theoretical digression of some interest.

\section{Analytical properties of the invariant measure}
\label{secanal}

The theory developed in the previous section can be also employed for an ambitious goal: to study numerically the analytical properties of the invariant measure $\mu$ of a $(\delta,\sigma)$--IFS. To give an idea of what we believe can be achieved along this line, in this section we briefly discuss two types of results, one for measures whose support is a full interval, the other for measured supported on Cantor sets.

The continuity properties of the measure $\mu$ follow in a complicated way from those of $\sigma$ and from the contraction ratio $\delta$. For instance, even in the realm of conventional, finite cardinality I.F.S., cases of point measures $\sigma$ leading to either singular continuous, or absolutely continuous measures $\mu$ are well known. Reference \cite{nalgo2} is an attempt to attack this problem in full generality: typically, we find that $\mu$ is ``more continuous'' than $\sigma$. We have proven that when $\sigma$ is absolutely continuous with a bounded density, so is $\mu$, for any $\delta$. This covers most cases commonly encountered in numerical analysis, but not the general situation that we discuss in this paper. Various techniques have been proposed in \cite{nalgo2} to verify numerically the continuity type of a measure $\mu$. We now present a different one.

Recall that the Nevai class of measures $N(a_\infty,b_\infty)$ contains the orthogonality measures associated with Jacobi matrices for which $a_n \rightarrow a_\infty$ and $b_n \rightarrow b_\infty$, as $n \rightarrow \infty$. From the speed of this convergence one can infer the continuity properties of $\mu$ (see e.g. \cite{vanas},\cite{damasi}). For instance, when
 \beq
 \sum_n |a_n-a_\infty| + | b_n - b_\infty | < \infty
 \label{nevai1}
 \end{equation}
the measure $\mu$ is absolutely continuous w.r.t. the Lebesgue measure on the interval $[a_\infty-2b_\infty,a_\infty+2b_\infty]$.

The numerical stability and the performance featured by Algorithm 2 permit to compute Jacobi matrices of large orders and therefore to hint to the existence of $a_\infty$ and $b_\infty$ and to the validity of eq. (\ref{nevai1}). We now apply this technique to the Erd\H{o}s problem of infinite Bernoulli convolutions. For details of this problem see the review paper \cite{berno} and references therein.
\begin{example}
\label{exa-berno}
Let $\sigma = \frac{1}{2} (D_{0}+D_{1})$ (a minor variation of example \ref{exa-mupp})
and select three values of $\delta$,  $\delta_1 = 2^{-1/2} \sim 0.7071067811865$, $\delta_2=3/4=.75$
and $\delta_3 =1/p_1 \sim 0.7548776662467$, $p_1$ being a Pisot number.
\end{example}

It is rigorously known that the generated measure $\mu$ is absolutely continuous in the first case and singular continuous in the third. It is also absolutely continuous for Lebesgue almost all values of $\delta$ between one half and one, but it is only conjectured that rational values in this interval, such as $\delta_2$, belong to this case. Notice that $\delta_2$ is very  close to $\delta_3$. Also notice that $b_\infty = 1/2$ is exactly known in the three cases.

\begin{figure}[tb]
\centerline{\includegraphics[width=6cm,height=12cm,angle=270]{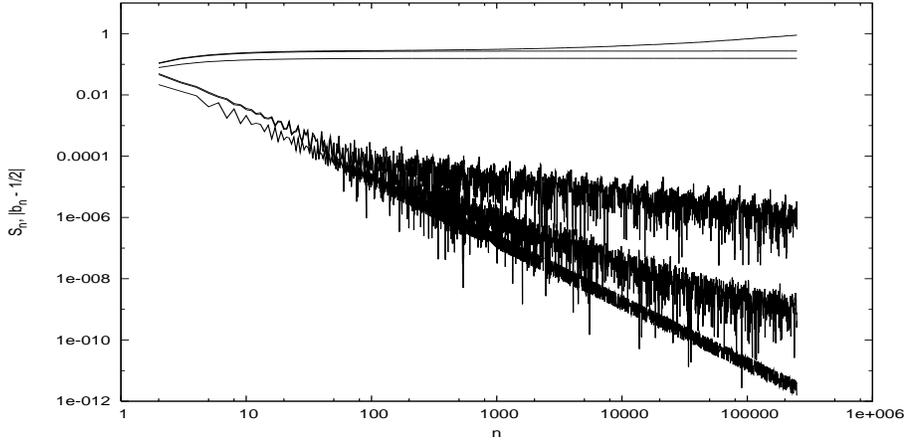}}
\caption{Differences $| b_n - \frac{1}{2} |$  (decreasing curves) and partial sums $S_n$ (increasing curves) versus $n$ for the three measures described in the text. Inside each group the curves arrange themselves for large $n$ from bottom to top, starting from Ex. \ref{exa-mupp} with $\delta=\delta_1$, followed by $\delta=\delta_2$ and finally $\delta=\delta_3$. } \label{nevai1b.fig}
\end{figure}

In figure \ref{nevai1b.fig} we plot both $| b_n - \frac{1}{2} |$ and the partial sums $S_n := \sum_j | b_j - 1/2 |$ versus $n$ in doubly logarithmic scale for the three cases listed. Convergence of $S_n$  is observed in the first two cases, as can be inferred by the exponent of the power-law decay of $| b_n - \frac{1}{2} |$. In the third case, this decay (which is also evident) is too slow to imply the absolute continuity of $\mu$. Therefore, our technique gives results that are consistent with rigorous facts and, which is more important, with the conjectured absolute continuity for the case $\delta = \delta_2$.

The measures just discussed are supported on the full interval $[0,1]$. Suppose now that $S_\mu$ is a Cantor set. The measure $\mu$ is then singular continuous and does not belong to a Nevai class. Our theory permits to derive estimates for the {\em capacity} of $S_\mu$, in the sense of potential theory \cite{stahl}. Precisely, we have

\begin{proposition}
Suppose that the distribution of fixed points
$\sigma$ and the invariant measure $\mu$ are regular, in the sense of potential theory. Let then $C_\sigma$
and $C_\mu$ be the capacities of the supports of $\sigma$ and of $\mu$, respectively. Then,
 \beq
  C_\sigma\leq   C_\mu \leq  C_\sigma  + \log (\bar{\delta}^{-1}) \leq C_\sigma  + \delta.
 \label{coreq1}
 \end{equation}
 \label{lemmacap}
\end{proposition}
{\em Proof.}
Consider the coefficient $h_n(\mu)$ of $s^n$ in $p_n(\mu;s)$: from eq. (\ref{nor2}) it follows that
 \beq
  h_n(\mu) = 1/ \prod_{j=1}^n b_j(\mu).
 \label{gaxx}
 \end{equation}
If $\mu$ is any regular probability measure, the
asymptotic relation $h_n(\mu) \sim e^{n C_\mu}$ holds \cite{stahl}. Then,
$\frac{h_n(\mu)}{h_n(\sigma)} $ behaves asymptotically as $e^{n
\Delta}$, where $\Delta = C_\mu - C_\sigma$. Since the support of
$\sigma$ is enclosed in that of $\mu$ this difference is always positive. Since $|\Omega^n_{k,r}| \leq 1$ for all $k,r$, using
Lemma \ref{lemgnorm}, we get
 \[
  h_n(\mu) \leq h_n(\sigma) \bar{\delta}^{-n},
  \]
and the result follows. \qed

\begin{remark}
It is possible to prove the regularity of I.F.S. singular continuous measures $\mu$ in large generality, see \cite{stahl}.
\end{remark}

Observe finally that
computing large order Jacobi matrices can lead to a numerical estimate of the capacity of the support of $\mu$, via the quantity $\frac{1}{n} \log (h_n(\mu)) = -\frac{1}{n} \sum_{j=1}^{n} \log(b_j(\mu))$ that converges to  $C_\mu$. Convergence is typically slow, but it can be accelerated by suitable extrapolation techniques \cite{claudemichela}.

\section{A spectral data technique}
\label{secspec}

The numerical determination of the Jacobi matrix of $\bar{\eta}$ in
Sect. \ref{sechomo}, eq. (\ref{bala2x}) can also be effected without
recurring to the algebraic theory developed in Sect. \ref{secphi} and leading to algorithms 1 and 2. In fact, the problem can be
cast into a forward/inverse Gaussian quadrature determination.
Nonetheless, this approach  only provides us with an analogue of
the recursive algorithm 1-Fix and not of the faster algorithm 2.
The key observation is the following:

\begin{lemma} \label{lemgau2}
The formula
 \beq
   \int d \bar{\eta}(s) f(s) =
   \sum_{j=1}^n  \sum_{k=1}^n  w_j^{(n)}(\eta) w_k^{(n)}(\sigma)
    f(\delta x_j^{(n)}(\eta) + \bar{\delta}  x_k^{(n)}(\sigma)),
 \label{gamma2}
 \end{equation}
where  $x_j^{(n)}(\cdot)$ and $w_j^{(n)}(\cdot)$ are Gaussian points and weights, respectively,
is exact for $f \in P_{2n-1}$. Therefore, it can be used to compute
$J_{\bar{\eta}}^{(n)}$ exactly.
\end{lemma}

{\em Proof.}
Let $f$ in eq. (\ref{bala2x}) be a
polynomial of degree at most $2n-1$ in the variable $s$. Then, $f
(\phi(\beta,s)) = f(\delta s + \bar{\delta} \beta)$ can be exactly
integrated with respect to the measure $\eta$ by $n$-points Gaussian
summation in the variable $s$. This Gaussian formula can be easily
obtained by the spectral problem of $J_{\eta}^{(n)}$. Also, $f
(\phi(\beta,s))$ is a polynomial of degree $2n-1$ in the variable
$\beta$, and the same is its integral w.r.t. $d\eta(s)$. This latter
polynomial can be exactly integrated with respect to the measure
$\sigma$ by an $n$-points Gaussian summation obtained from
$J_{\sigma}^{(n)}$. \qed

The above Lemma can also serve as an alternative proof of Thm.
\ref{theodirect2}. Next, observe that the r.h.s. of eq. (\ref{gamma2}) is the integral of $f$ with respect to the sum of $n^2$ atomic measures. Stable algorithms for computing the Jacobi matrix of a finite sum of atomic measures, due among others to De Boor and Golub \cite{gene1}, Gragg and Harrod \cite{gragg}, Fischer \cite{hans3}, Reichel \cite{lotrei} and Laurie \cite{dirk9} are well known and can be put to use, to give the following:
\begin{itemize}
\item[] {{\bf Algorithm 3. Computing the I.F.S. convolution.} \\
 {\bf Input}: the (truncated) Jacobi matrices $J_\sigma^{(\bar{n})}$  and $J^{(\bar{n})}_{{\eta}}$, the contraction factor $\delta$, the truncation size $\bar{n}$.
 \\ {\bf Output}: the (truncated) Jacobi matrix of  $\bar{\eta} := \Phi_\delta(\sigma;\eta)$.}
\item[1:] Compute Gaussian points and weights for  $\sigma$
from  $J_{\sigma}^{(\bar{n})}$.
\item[2:] Compute Gaussian points and weights for $\eta$
from $J_{\eta}^{(\bar{n})}$.
\item[3:] Using Lemma \ref{lemgau2} compute $J_{\bar{\eta}}^{(\bar{n})}$ using one of the algorithms just quoted.
\end{itemize}

It is immediate to obtain an iterative version, 3-Fix, along the same lines of Alg. 1-Fix: from step [3] loop back to step [2] replacing $\eta$ by $\bar{\eta}$. This fixed--point algorithm then provides us with the Jacobi matrix of $\mu$.
As a matter of facts, algorithm 3-Fix works fine as far as absolutely continuous measures $\sigma$ (like that in Example \ref{exa-muac}) are involved. Instead, when $\mu$ ({\em not} $\sigma$) is supported on a Cantor set (like {\em e.g.} in Example \ref{exa-mupp}), we have observed that it achieves convergence only for Jacobi matrices of the size of about a thousand. This is due to the fact that the $n^2$ points in Lemma \ref{lemgau2} crowd around a fractal and the relative precision in their distance diminishes. As explained in detail by Laurie \cite{dirk9} (see also \cite{gragg,leary}) this fact impairs the reconstruction of the Jacobi matrix from the Gaussian points and weights. One can therefore appreciate by comparison the computational advantage brought about by the algebraic theory of Sect. \ref{secphi}.

\section{I.F.S. Quadratures and the Inverse Problem}
\label{secinve}

The algorithm 2 presented above can be reversed, in order to compute $J_\sigma$ from $J_\mu$. This is the basis of the solution of an inverse problem, that can be used in an approximation problem: that of finding I.F.S. quadratures \cite{steve,nalgo1}.

\begin{definition}  \label{ggaus1}
Given a target measure $\mu$, whose support is enclosed in a finite interval, an I.F.S. quadrature for $\mu$ is a sequence of I.F.S. measures $\mu^{(n)}$ that satisfies
 $
    J^{(n)}_{\mu^{(n)}} = J^{(n)}_\mu,
$
for any $n \in {\bf N}$.
\end{definition}

\begin{remark}
Def. (\ref{ggaus1}) implies that $\mu^{(n)}$ integrates exactly polynomials up to degree $2n-1$ and therefore the sequence $\{\mu^{(n)}\}$ is weakly convergent to $\mu$. Clearly, the linear combination, with Gaussian weights, of the atomic measures sitting at the Gaussian points of order $n$ is an I.F.S. quadrature, degenerate in the sense that Gaussian points are the fixed points of a finite set of maps with contraction rate $\delta=0$.
\end{remark}

Def. (\ref{ggaus1}) formalizes a {\em truncated inverse problem}, in the family of fractal inverse problems \cite{ba0,ba2,steve,gio1,ed}. If we now restrict ourselves to $(\delta,\sigma)$--I.F.S. of the kind (\ref{mappi}), the following approximation result guarantees that solutions do exist:
 \begin{theorem}[\cite{gio1}]
Let $\mu$ be a measure with an infinite number of points of
increase. Then, for all $n > 0$ there exists  $\delta_n(\mu) > 0$ such
that for all $\delta \in [0,\delta_n(\mu))$ there exists an homogeneous
affine I.F.S. with $n$ maps that satisfies Def. (\ref{ggaus1}).
 \label{thgio}
  \end{theorem}

This theorem and Def. \ref{ggaus1} also imply that {\em any} finite symmetric tridiagonal matrix with positive out--diagonals $b_n$ is the truncation of the Jacobi matrix of a $(\delta,\sigma)$--I.F.S. with non-vanishing $\delta$. We shall now apply this theorem while developing the inverse of
algorithm 2, in a form that is also suitable for the numerical determination of the maximal value $\delta_n(\mu)$.

\begin{theorem}[\cite{nalgo1}] The truncated Jacobi matrix $J^n_\sigma$ of the distribution of fixed points $\sigma$ of a
$(\delta,\sigma)$--I.F.S. with contraction ratio $\delta$
that provides an I.F.S. quadrature of a measure $\mu$
can be computed recursively from the
truncated Jacobi matrix $J^n_\mu$, provided $\delta \leq \delta_n(\mu)$. Whether the last condition holds can be verified recursively. \label{theoinverse}
\end{theorem}
{\em Proof.}
The algorithm is effected in the following sequence of steps
\begin{itemize}
\item[] {{\bf Algorithm 4. Solving the inverse I.F.S. problem.} \\
 {\bf Input}: the (truncated) Jacobi matrix $J_\mu^{(\bar{n})}$  of the target measure $\mu$, the contraction ratio $\delta$, the maximum size $\bar{n}$.
 \\ {\bf Output}: the (truncated) Jacobi matrix of $\sigma$, the largest allowed truncation size $\hat{n}$ }.
\item[0:] Initialization: $n=0$. One has $\Omega^0_{0,0}= 1$, since
$p_0(\mu;s)=p_0(\sigma;\beta)=1$, and $b_0(\mu)=b_0(\sigma)=0$.
\item[1:] Induction hypothesis:
$\{\Omega^j,j=0,\ldots,n\}$, $\{a_j(\sigma), j=0,\ldots,n-1 \}$, and
$\{b_j(\sigma), j=0,\ldots,n \}$ are known.
\item [2:] Computation of $a_n(\sigma)$: Lemma \ref{lemgamm3}.
Observe that $a_n(\sigma)$ has a non-zero coefficient in eq. (\ref{ad2d}), due to Lemma \ref{lg1e}.
\item [3:] Computation of the matrix $\tilde{\Omega}^{n+1}$: Lemma \ref{lgammarec}.
    \item [4:] Computation of $b^2_{n+1}(\sigma)$: Lemma \ref{lemrnp2x}, eq. (\ref{gam356}).
\item [4:] Stopping alternative: either $b^2_{n+1}(\sigma) > 0$, therefore continue,
or else
 $\delta > \delta_{n+1}(\mu)$, $\hat{n} = n$ and stop.
\item [6:] Computation of $\Omega^{n+1}$: divide
$\tilde{\Omega}^{n+1}$ by $b_{n+1}(\mu)$.
\item [7:] If $n<\bar{n}$ augment $n$ to $n+1$ and loop back to 1, else $\hat{n} = \bar{n}$ and stop.
\end{itemize}

\begin{remark}
When termination occurs at a certain value of $n = \hat{n} < \bar{n}$ at step 4, then $\delta$ is larger than $\delta_{n+1}(\mu)$, but smaller than $\delta_{n}(\mu)$. Therefore, using Algorithm 4 iteratively at different values of $\delta$, one can determine the sequence of values $\delta_n$, at varying $n$. In principle, the algorithm never stops only if the target measure $\mu$ is {\em exactly} generated by an affine IFS with contraction ratio $\delta$ and a measure $\sigma$ with an infinite number of points of increase.
\end{remark}

To establish the numerical stability of both the forward algorithm 2 and the reverse algorithm 4 we have chosen a target Jacobi matrix of particular significance, the Fibonacci tridiagonal matrix, whose orthogonality measure is singular continuous. This measure is {\em not} the invariant measure of a $\delta$-homogeneous I.F.S. and yet, as seen above, any finite truncation of its Jacobi matrix coincides with the truncation of the Jacobi matrix of a $(\delta,\sigma)$--I.F.S.

\begin{example}
Let $a_n = 0$ for all $n$, and let $b_n$ take either the value $A=2/5$ or the value $B=1/2$. These values are arranged in the aperiodic Fibonacci sequence $ABABBA \dots$, generated by the substitution rules $A \rightarrow AB$, $B \rightarrow A$ on the seed $A$. It has the property that two A's never follow each other, but are separated by at most two B's.
\label{exa-fibo}
\end{example}

\begin{figure}[tb]
\centerline{\includegraphics[width=6cm,height=12cm,angle=270]{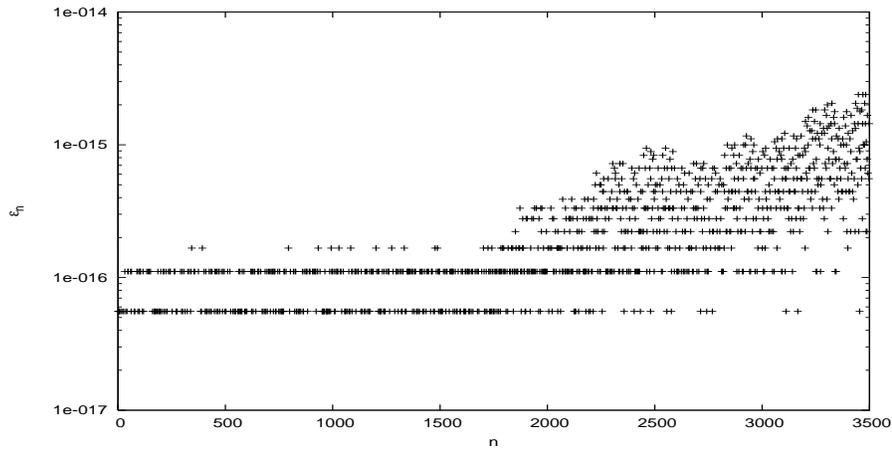}}
\caption{Example \ref{exa-fibo}. Absolute errors $\varepsilon_n$ versus $n$, for the reconstruction of the out-diagonal Fibonacci Jacobi matrix coefficients $b_n$, at $\delta = 1.119837 \; 10^{-6}$, the maximum allowed value for $\delta$ at ${\bar{n}}=3500$ being computed as approximately $\delta_{3500}(\mu) = 1.124611 \; 10^{-6}$.}\label{diffefib.fig}
\end{figure}

Using algorithm 4 we have computed the sequence $\delta_n(\mu)$ at increasing values of $n$ up to $\bar{n}$, as well as the Jacobi matrix $J^{(\bar{n})}_\sigma$ for a feasible value of $\delta$, close to the maximum allowed value $\delta_{\bar{n}}(\mu)$. We have then applied Algorithm 2 to recompute the original Fibonacci Jacobi matrix. While the null diagonal entries are recovered exactly because of the nature of the algorithms, Figure \ref{diffefib.fig} plots the absolute errors $\varepsilon_n$ in the reconstruction of the sequence of $b_n$. The observed behavior, confirmed by other experiments, is the most convincing experimental verification of the numerical stability of the direct and inverse algorithms 2 and 4 presented in this paper.

\section{Conclusions}
\label{secconclu}

We have presented a new family of algorithms for the direct/inverse computation of the Jacobi matrix of the invariant measure of a homogeneous affine I.F.S.. Experimental results suggest that these algorithms are stable to large orders when programmed in floating point arithmetics.

On the one hand, this remarkable stability calls for a detailed error analysis, that should unveil the reasons why the algebraic treatment of Sect. \ref{secphi} is more stable than any other existing technique (to the author knowledge) and in particular than the Gaussian technique of Sect. \ref{secspec}. We conjecture that this is due to the fact that in our approach we exclusively deal with Jacobi matrices, but a more thorough investigation, that we plan to develop in further publications, is in order.

On the other hand, the versatile tools that have been introduced in this paper can now be applied to a variety of problems, both from the theoretical and from the applied side. In the first respect, we would like to investigate to what extent we can infer the fine structure properties of the generated measure $\mu$ from those of $\sigma$ and whether more can be said from the potential theoretical point of view, {\em e.g.} on the asymptotic properties of the sequence of Jacobi matrix entries and on the Fourier transform of $\mu$ and of its orthogonal polynomials. In the second respect, we would like to evaluate the full potential of the approximation/inverse problem of Sect. \ref{secinve} on significant problems, until now beyond the reach of conventional algorithms.

{\em Acknowledgements}
It is a pleasure to have this opportunity to thank Dirk Laurie for providing his code for \cite{dirk9} and for related interesting discussions.

\end{document}